\documentclass[twoside,11pt]{article}

\usepackage{jmlr2e}
\usepackage{bbm}
\usepackage{amsfonts}
\usepackage{amsmath}
\usepackage{bm}
\usepackage{hyperref}
\hypersetup{hypertex=true,
            colorlinks=true,
            linkcolor=blue,
            anchorcolor=blue,
            citecolor=blue}

\jmlrheading{}{}{}{}{}{Xiao Fang and Malay Ghosh}

\ShortHeadings{Bernstein von-Mises Theorem for $g$-Priors and Nonlocal
Prior}{Fang  and Ghosh}
\firstpageno{1}

\begin{document} 
\title{Bernstein von-Mises Theorem for g-prior and nonlocal
prior}

\author{\name Xiao Fang \email xiaofang@ufl.edu\\
       \name Malay Ghosh \email ghoshm@ufl.edu \\
       \addr Department of Statistics\\
       University of Florida\\
       Gainesville, FL 32611, USA}

\editor{}

\maketitle

\begin{abstract}
The paper develops Bernstein
von Mises Theorem under hierarchical g -priors for linear regression models.
The results are obtained both
when the error variance is known,
and also when it is unknown. An
inverse gamma prior is attached to
the error variance in the later case.
The paper also demonstrates some
connection between the total
variation and $\alpha$-divergence measures.
\end{abstract}

\section{ Introduction}
Bernstein von-Mises(BVM) Theorem occupies a central place in the statistics literature. In its original version, the theorem establishes asymptotic normality of the posterior distribution of the parameter vector of interest for a regular finite dimensional model.The asymptotic distribution is centered at the posterior mean or posterior mode with inverse of the observed information matrix as its variance-covariance matrix.\\

There exists a huge literature establishing BVM theorems under  different scenarios, but still under the framework of a regular finite dimensional model. It is only recently that attention has bean shifted to establishing  BVM theorems where the number of covariates can even depend on the sample size. This has been demonstrated by \citep{Song:2017} for a linear regression model with a scale-mixed normal prior with dominating peak around zero. The theorem has been proved with spike and slab prior by \citep{Castillo:2015}, and \citep{Martin:2020} under a linear regression model.\\

The objective of our paper is to establish BVM theorems both under the g-prior\citep{Zellner:1986,Liang:2008},  and non-local priors \citep{Johnson:2012} once again in the set up of linear regression model. For g-priors we have proved the result both when the error variance is known or unknown. For non-local prior, we have established the result when the error is known.\\

The BVM theorems are usually established by either showing that the total variation between the normalized posterior density and the normal density converges to zero as the sample sizes tend to infinity or the Hellinger distance  between the two densities  converges to  zero as the sample size tend to infinity. Indeed the two  criteria are  equivalent as we will demonstrate in the next section. We have adopted  in this paper the Hellinger distance criterion for establishing the BVM Theorem.\\

The outline of the remaining  sections is as follows. In Section 2, we have proved the equivalence of density convergence  under the Renyi divergence and the total variation  divergence criteria. As we shall see, the Hellinger divergence is  a special case of the Renyi divergence. We have also stated  a  few lemmas related to tail tail bounds of chi-squared distribution, used repeatedly in the proofs of our main results for g-priors in Section 3 and for non-local priors in Section 4. Some final remarks are made in Section 5. The proofs all the results are given in Section 6.

\section{Some Preliminary Result}
\subsection{$\alpha$-divergence Criterion}
We first introduce the \citep{Renyi:1961} divergence criterion, more popularly known as the $\alpha$-divergence criterion between two densities $p_1$ and $p_2$, each with respect to some $\sigma$-finite measure $\mu$ as
\begin{displaymath}
D_{\alpha}(p_1,p_2)=[1-\int p_1^{\alpha}p_2^{1-\alpha} d \mu]/[\alpha(1-\alpha)],0<\alpha<1.
\end{displaymath}
It maybe noted that when $\alpha=\frac{1}{2}$,
\begin{displaymath}
D_{1/2}(p_1,p_2)=4[1-\int p_1^{1/2} p_2^{1/2} d \mu]=2\int (p_1^{1/2}-p_2^{1/2})^2d \mu=2H^2(p_1,p_2),
\end{displaymath}
where $H(p_1,p_2)=[\int (p_1^{1/2}-p_2^{1/2})^2 d\mu]^{1/2}$ is the Hellinger distance between two densities $p_1$ and $p_2$. As an aside, the two limitting case $\alpha \to 0$ and $\alpha \to 1$ yield respectively the Kullback-Leibler(KL) divergence $KL(p_2,p_1)=\int \log(p_2 /p_1) p_2d\mu$ and $KL(p_1,p_2)=\int \log(p_1 /p_2) p_1d\mu$. We also recall that the total-variation between two densities  $p_1$ and $p_2$  is given by 
$TV(p_1,p_2)=\frac{1}{2}\int |p_1-p_2| d \mu$.

\begin{lemma}
Consider a sequence of densities $p_n(n\ge1))$ and a density $p$ all with respect to some $\sigma$-finite  measure $\mu$. Then $TV(p_n,p)\to 0 \Rightarrow D_{\alpha}(p_n,p)\to 0$  for all $0<\alpha<1$.\label{lemma0}
\end{lemma}
Proof:  For $0<\alpha<1$,
\begin{displaymath}
\begin{split}
D_{\alpha}(p_n,p)&=[[1-\int p_n^{\alpha}p^{1-\alpha} d \mu]/[\alpha(1-\alpha)]]\\
& \le [1-\int \{\min(p_n,p)\}^{\alpha}\{\min(p_n,p)\}^{1-\alpha} d \mu]/[\alpha(1-\alpha)]\\
& = [1-\int \min(p_n,p) d \mu]/[\alpha(1-\alpha)]\\
&= [2-2\int \min(p_n,p) d \mu]/[2\alpha(1-\alpha)]\\
&= [p_n + p-2 \int \min(p_n,p) d \mu]/[2\alpha(1-\alpha)]\\
&=\frac{1}{2} [\int |p_n-p|d\mu] / [\alpha(1-\alpha)]=TV(p_n,p)/ [\alpha(1-\alpha)].
\end{split}
\end{displaymath}
Hence, $TV(p_n,p) \to 0$ implies $D_{\alpha}(p_n,p) \to 0$ for all $0<\alpha<1$.

\begin{remark}
The result is not necessarily true when $\alpha \to 0$ or $\alpha \to 1$.
\end{remark}

In view of the above result, $D_{1/2}(p_n,p)=2H^2(p_n,p)\to 0$ when $TV(p_n,p) \to 0$. But we also know that $TV(p_n,p) \le H(p_n,p)$.  Thus $H(p_n,p) \to 0 \Leftrightarrow TV(p_n,p) \Leftrightarrow D_{\alpha}(p_n,p)$ for all $0<\alpha<1$.

\subsection{Tail bounds for the Chisquared distribution}
In this subsection, we state a few  lemmas related to tail behavior of Chisquared distribution.

\begin{lemma}
For any $a>0$, we have
\begin{displaymath}
P(\vert\chi^2_p -p\vert>a) \le 2 \exp(-\frac{a^2}{4p}).
\end{displaymath} \label{lemma1}
\end{lemma}
\noindent
Proof: See Lemma 4.1 of \citep{Cao:2020}.

\begin{lemma} \noindent \\
(i) For any $c>0$, we have
\begin{displaymath}
P(\chi^2_p(\lambda)-(p+\lambda)>c) \le \exp(-\frac{p}{2} \{\frac{c}{p+\lambda}-log(1+\frac{c}{p+\lambda})\}).
\end{displaymath}
(ii) for $\omega <1$ then
\begin{displaymath}
P(\chi^2_p(\lambda) \le \omega \lambda) \le c_1 \lambda^{-1} \exp \{-\lambda(1-\omega)^2/8\},
\end{displaymath}
where $c_1>0$ is a constant.\\
\\
(iii)
For any $c>0$, $P(\chi^2_p(\lambda)-p \le -c) \le  \exp(-\frac{c^2}{4p})$.\label{lemma2}
\end{lemma} 
\noindent
Proof: See  Lemma 4.2 of \citep{Cao:2020}  for the proof of (i). (ii) is proved in \citep{Shin:2019}.  Finally, (iii)  follows from  the fact  if $U \sim \chi^2_p(\lambda)$, then $P(U>c)$ is strictly increasing in $\lambda$ for fixed $p$ and $c>0$. Hence
\begin{displaymath}
P(\chi^2_p(\lambda)-p \le -c) \le P(\chi^2_p-p \le -c) \le \exp(-\frac{c^2}{4p}).
\end{displaymath}

\subsection{Notations}
Let $\hat{\beta}=(X^T X)^{-1} X^T Y$. $H$ denotes the Hellinger distance. It's well known that $H^2$ a is convex function. Also, let $P_X=X(X^T X)^{-1} X^T$, the projection matrix. 

\section{Bernstein von-Mises theorem for g-prior}
Consider the standard linear regression model 
\begin{equation}
Y=X\beta+ \sigma z,  \label{linear model}
\end{equation}
where $Y$ is a $(n \times 1)$ vector of response variables, $X$ is $n \times p$ design matrix, $\beta$ is a $(p \times 1)$ regression parameter, $\sigma >0 $ is an  scale parameter and $z \sim N_n(0,I_n)$, $I_n$ being the identity matrix. We will allow $p\equiv p_n$ to depend on $n$\\

Following \citep{Liang:2008}, we consider the following hierarchical model 
\begin{equation}
\begin{split}
Y|\beta & \sim  N(X \beta,\sigma^2 I_n ) \\
\beta|g & \sim  N(0,g \sigma^2(X^T X)^{-1}) \\  
\pi(g ) & \propto  g^{b-1}(1+g)^{-(a+b)}, g>0. 
\end{split}
\end{equation}

Writing $\omega=({1+g})^{-1}$,  $\omega \sim Beta(a,b)$. Then we have
\begin{equation}
\pi(\omega|\sigma^2,Y) \sim \omega^{a+\frac{p}{2}-1}(1- \omega)^{b-1} \exp (-\omega \frac{Y^T P_X Y}{\sigma^2}),
\end{equation}

\begin{equation}
\pi(\beta|\omega, \sigma^2,Y) \sim \exp(-\frac{(\beta-(1-\omega)\hat{\beta})^T(X^T X)(\beta-(1-\omega)\hat{\beta})}{2(1-\omega)\sigma^2}).
\end{equation}
We will establish BVM theorem both when  $\sigma^2$ is fixed and when it is unknown .

\subsection{Fixed $\sigma^2$ case}
We can show that for fixed $\sigma^2$ case, the posterior distribution of $\omega$ will concentrate near $0$, so that  the posterior distribution of $\beta$ will tend to a normal distribution centered at $\hat{\beta}$ with variance matrix $\sigma^2(X^T X)^{-1}$.\\
\noindent 
\begin{lemma} 
Consider the  hierarchical prior given before  and assume $\sigma^2$ is fixed. Also we denote the true model as $Y \sim  N(X \beta_0,\sigma^2 I_n )$. Let  ${\beta_0}^T {X^T X} \beta_0 \asymp n$  and  $p=o(\frac{n}{\log n})$, $p \to \infty$ as $n \to \infty$.  For any diverging sequence $\xi_n$, let $t_n=\xi_n p \log n$. Then we have 
\begin{displaymath}
E_0[P(\omega > \frac{t_n}{n}|Y)]=o(1),
\end{displaymath}
where $E_0$ denotes expectation under $\beta =\beta_0$.
\label{lemma3}
\end{lemma}
\noindent
By applying Lemma \ref{lemma3}, we get BVM theorem for g-prior with fixed $\sigma^2$.

\begin{theorem}
Let $p=O(n^{\frac{1}{4}}$) and $\xi_n= \log n$ in Lemma \ref{lemma3}. Then  under  conditions of Lemma \ref{lemma3}, 
\begin{displaymath}
E_0 H^2(\pi(\beta|\sigma^2,Y),N(\hat{\beta}, \sigma^2(X^T X)^{-1})) \to 0 ,
\end{displaymath}
as $n \to \infty$.  \label{theorem1}
\end{theorem}

\subsection{Unknown $\sigma^2$ case}
Now, we assume inverse gamma prior  $
\pi(\sigma^2) \sim IG(\frac{c}{2},\frac{d}{2})$. 
\\
 
Then we have 
\begin{equation}
\beta| \omega ,Y \sim t_{n+c}((1-\omega)\hat{\beta},\frac{d+Y^T(I-P_X) Y+\omega Y^T P_X Y}{(1-\omega)(n+c)}(X^TX)^{-1}),
\end{equation}
where $t_{n+c}$  is a student's t-distribution with $n+c$ degrees of freedom,
\begin{equation}
\pi(\omega|Y) \propto (d+Y^T(I-P_X) Y+\omega Y^T P_X Y)^{-\frac{n+c}{2}}\omega^{a+\frac{p}{2}-1}(1- \omega)^{b-1}.
\end{equation}
\noindent
\\
Similar to Lemma \ref{lemma3}, we can also show the posterior distribution of $\omega$ will concentrate near $0$.
\begin{lemma} \noindent \\
Consider the  hierarchical prior given before  and assume $\pi(\sigma^2) \sim IG(\frac{c}{2},\frac{d}{2})$. Let ${\beta_0}^T {X^T X} \beta_0 \asymp n$  and  $p=o(\frac{n}{\log n})$, $p \to \infty$ as $n \to \infty$.  Then for any diverging sequence $\xi_n$, let $t_n=\xi_n p \log n$, we have 
\begin{displaymath}
E_0[P(\omega > \frac{t_n}{n}|Y)]=o(1).
\end{displaymath}
\label{lemma4}
\end{lemma}
\noindent
By applying Lemma \ref{lemma4}, we get  BVM theorem for the  g-prior with unknown $\sigma^2$.

\begin{theorem}
Under the same conditions of Lemma \ref{lemma4},  let $p=O(n^{\frac{1}{4}}$) and $\xi_n= \log n$. Then we have 
\begin{displaymath}
E_0 H^2(\pi(\beta|\sigma^2,Y),N(\hat{\beta}, \sigma_0^2(X^T X)^{-1})) \to 0 ,
\end{displaymath}
as $n \to \infty$.  \label{theorem2}
\end{theorem}

\section{Bernstein von-Mises theorem for nonlocal priors}
Now we consider nonlocal prior in \citep{Johnson:2012}, the prior on $\beta$ is a pMoM density, as considered in their paper specifically,  
\begin{displaymath}
\pi(\beta|\tau, \sigma^2, r ) =d_p (2 \pi)^{-\frac{p}{2}}(\tau\sigma^2)^{-rp-\frac{p}{2}}|A_p|^{1/2}\exp[-\frac{1}{2 \tau \sigma^2} \beta^T A_p \beta] \prod_{i=1}^p \beta_i^{2r},
\end{displaymath}
for $\tau >0$, $A_p$  a  $p \times p$ nonsingular scale matrix,  $r=1,2, \cdots$, and $d_p$ is a normalizing constant.\\
Then we have 
\begin{displaymath}
\pi(\beta_k| k,Y,\tau, r, \sigma^2)  
= E_k^{-1}(\Pi_{i=1}^{k} \beta_i^{2r})(2 \pi \sigma^2)^{-\frac{k}{2}}|C_k|^{-\frac{1}{2}}(\Pi_{i=1}^{k} \beta_i^{2r}) \exp[\frac{(\beta_k-\tilde{\beta}_k)^{T}C_k(\beta_k-\tilde{\beta}_k)}{2 \sigma^2}],
\end{displaymath}
where $C_k=X_k^T X_k +\frac{1}{\tau} A_k$, $\tilde{\beta}_k=C_k^{-1}X_k^{T}Y$, for $k=\{k_1,k_2, \cdots, k_j\}(1 \le k_1 < \cdots < k_j \le p)$ , we let $X_k$ denote the design matrix formed from
the columns of $X_n$ corresponding to model $k$,  $E_k(\cdot)$ denotes expectation with respect to multinormal distribution with mean $\tilde{\beta}_k$ and covariance matrix $\sigma^2 C_k^{-1}$.\\
\\
We have the following   BVM theorem for $\sigma^2$ is fixed case.
\begin{theorem}
Under the conditions of Theorem 1 in \citep{Johnson:2012},  we have 
\begin{displaymath}
E_{\beta_t^0} H^2(\pi(\beta|\tau, r, \sigma^2, Y),N(\hat{\beta}_t, \sigma^2(X_t^T X_t)^{-1})) \to 0 ,
\end{displaymath}
as $n \to \infty$.  \label{theorem3}
\end{theorem}

\section{Discussion}
In this paper, we establish BVM theorems both under g-prior \citep{Zellner:1986,Liang:2008} and non-local priors(\citep{Johnson:2012}) once again in the set up of linear regression model. For the g-prior we have proved the result both when the error variance is known or unknown. For non-local prior, we have established the result when the error is known. However for non-local priors, BVM theorem for  unknown error variance   remains to be established.\\
\\
With   BVM theorem, we are also able to obtain  many interesting results like valid uncertainty quantification as in \citep{Martin:2020}, and we will discuss it in future paper.

\section{Proofs }

When $\beta=\beta_0$
\begin{displaymath}
\frac{Y^T P_X Y}{\sigma^2} \sim {\chi}^2_p(\frac{{\beta_0}^T {X^T X} \beta_0 }{2 \sigma^2}).
\end{displaymath}
Since  ${\beta_0}^T {X^T X} \beta_0 \asymp n$, $p=o(n)$, $p \to \infty$  as $n \to \infty$, by Lemma \ref{lemma2} we have 
\begin{displaymath}
P(\frac{1}{2} \frac{{\beta_0}^T {X^T X} \beta_0 }{2 \sigma^2}   \le \frac{Y^T P_X Y}{\sigma^2}   \le p+2 \frac{{\beta_0}^T {X^T X} \beta_0 }{2 \sigma^2}) \to 1
\end{displaymath}
as $n \to \infty$, also there exist constant $l_1, l_2>0$, such that 
\begin{displaymath}
P(l_1 n  \le \frac{Y^T P_X Y}{\sigma^2}   \le l_2 n ) \to 1
\end{displaymath}
 as $n \to \infty$.\\
\\
Define  $S_n=\{Y: l_1 n  \le \frac{Y^T P_X Y}{\sigma^2}   \le l_2 n\}$, then $P(S_n) \to 1$, as $n \to \infty$.\\
\\
\textbf{Proof of Lemma \ref{lemma3}:}\\
Throughout the proof, we will use the notation $C(>0)$ as a generic constant.\\

We have $P(\omega > \frac{t_n}{n}|Y) \le \textbf{1}_{S_n}P(\omega > \frac{t_n}{n}|Y)+ \textbf{1}_{S^c_n}$ and $E_0\textbf{1}_{S^c_n}=P(S_n^c) \to 0$ as $n \to \infty$. So it suffices to show  that
\begin{displaymath}
E_0\{\textbf{1}_{S_n}P(\omega > \frac{t_n}{n}|Y)\} \to 0,
\end{displaymath}
as $n \to \infty$. In view of $\int_0^{1}\pi(\omega|Y)d \omega  \ge \int^1_{\frac{t_n}{n}}\pi(\omega|Y)d \omega$, all we need to show is
\begin{displaymath}
\frac{\int_0^{\frac{t_n}{n}}\pi(\omega|Y)d \omega}{\int^1_{\frac{t_n}{n}}\pi(\omega|Y)d \omega} \textbf{1}_{S_n}  \to \infty  ,
\end{displaymath}
as $n \to \infty$.\\

We have \begin{displaymath}
\pi(\omega|Y)=\frac{1}{m(Y)} \omega^{a+\frac{p}{2}-1}(1- \omega)^{b-1} \exp (-\omega \frac{Y^T P_X Y}{\sigma^2}),
\end{displaymath}
 where $m(Y)= \int_0^1 \omega^{a+\frac{p}{2}-1}(1- \omega)^{b-1} \exp (-\omega \frac{Y^T P_X Y}{\sigma^2}) d \omega$.
\begin{equation}
\begin{split}
& \textbf{1}_{S_n}{\int_0^{\frac{t_n}{n}}\omega^{a+\frac{p}{2}-1}(1- \omega)^{b-1} \exp (-\omega \frac{Y^T P_X Y}{\sigma^2})d \omega} \\
\ge & {\int_0^{\frac{t_n}{n}}\omega^{a+\frac{p}{2}-1}(1- \omega)^{b-1} \exp (-l_2 n \omega)d \omega}\\
\ge & C {\int_0^{\frac{t_n}{n}}\omega^{a+\frac{p}{2}-1} \exp (-l_2 n \omega)d \omega}\\
\ge & C {\int_0^{\frac{1}{n}}\omega^{a+\frac{p}{2}-1} \exp (-l_2 n \omega)d \omega}\\
\ge & C  {\int_0^{\frac{1}{n}}\omega^{a+\frac{p}{2}-1} d \omega}=\frac{C^{\prime \prime}}{a+\frac{p}{2}}n^{-(a+\frac{p}{2})}.  
\end{split} \label{numerator1}
\end{equation}
We also have
\begin{equation}
\begin{split}
&\textbf{1}_{S_n}{\int^1_{\frac{t_n}{n}}\omega^{a+\frac{p}{2}-1}(1- \omega)^{b-1} \exp (-\omega \frac{Y^T P_X Y}{\sigma^2})d \omega} \\
\le & {\int^1_{\frac{t_n}{n}}\omega^{a+\frac{p}{2}-1}(1- \omega)^{b-1} \exp (-\omega l_1 n)d \omega} \\
\le &  \exp (- l_1 t_n) {\int^1_{0}\omega^{a+\frac{p}{2}-1}(1- \omega)^{b-1} d \omega}\\
\le & \frac{\Gamma(a+\frac{p}{2})\Gamma(b)}{\Gamma(a+\frac{p}{2}+b)} \exp (- l_1 t_n) \le \exp (- l_1 t_n).
\end{split}  \label{denom 1}
\end{equation}
Combining (\ref{numerator1}) and (\ref{denom 1}), we obtain 
\begin{displaymath}
\frac{\int_0^{\frac{t_n}{n}}\pi(\omega|Y)d \omega}{\int^1_{\frac{t_n}{n}}\pi(\omega|Y)d \omega} \textbf{1}_{S_n} \ge \frac{C^{\prime \prime}}{a+\frac{p}{2}} \frac{\exp (- l_1 t_n)}{n^{-(a+\frac{p}{2})}} \to \infty.
\end{displaymath}
This proves Lemma 3. \\
\\
\\
\textbf{Proof of Theorem \ref{theorem1}:}\\
 $\pi(\beta|\sigma^2,Y)=\int_0^1\pi(\beta|\omega,\sigma^2, Y)\pi(\omega|\sigma^2,Y)d\omega$.\\
By convexity of $H^2$, Jensen inequality and $\beta|\omega,\sigma^2, Y \sim N((1-\omega)\hat{\beta},(1-\omega)\sigma^2(X^T X)^{-1})$, we have 
\begin{displaymath}
\begin{split}
& H^2(\pi(\beta|\sigma^2,Y),N(\hat{\beta}, \sigma^2(X^T X)^{-1}))\\
\le & \int_0^1\pi(\omega|\sigma^2,Y) H^2(N((1-\omega)\hat{\beta},(1-\omega)\sigma^2(X^T X)^{-1}),N(\hat{\beta},\sigma^2(X^T X)^{-1}))d \omega\\
\le &\int_0^{\frac{t_n}{n}}\pi(\omega|\sigma^2,Y) H^2(N((1-\omega)\hat{\beta},(1-\omega)\sigma^2(X^T X)^{-1}),N(\hat{\beta},\sigma^2(X^T X)^{-1}))d \omega+\int_{\frac{t_n}{n}}^1\pi(\omega|\sigma^2,Y) d \omega\\
= &\int_0^{\frac{t_n}{n}}\pi(\omega|\sigma^2,Y) H^2(N((1-\omega)\hat{\beta},(1-\omega)\sigma^2(X^T X)^{-1}),N(\hat{\beta},\sigma^2(X^T X)^{-1}))d \omega+P(\omega>\frac{t_n}{n}|Y)\\
= &I+II(say).
\end{split}
\end{displaymath}
By Lemma \ref{lemma3}, we have $E_0 II \to 0$. So it suffices to show  that $E_0 I \to 0$. \\

By \citep{Pardo:2006,Ghosh:2008} and , 
\begin{displaymath}
\begin{split}
& H^2(N((1-\omega)\hat{\beta},(1-\omega)\sigma^2(X^T X)^{-1}),N(\hat{\beta},\sigma^2(X^T X)^{-1}))\\
= &1- \frac{(1-\omega)^{\frac{p}{4}}}{(1-\omega/2)^{\frac{p}{2}}} \exp\{ -  \frac{\omega^2}{4(2-\omega)} \frac{Y^TP_X Y}{\sigma^2}\}.
\end{split}
\end{displaymath}
Thus 
\begin{displaymath}
\begin{split}
& I= \int_0^{\frac{t_n}{n}}\pi(\omega|\sigma^2,Y) [1- \frac{(1-\omega)^{\frac{p}{4}}}{(1-\omega/2)^{\frac{p}{2}}} \exp\{ - \frac{1}{4}\cdot  \frac{\omega^2}{2-\omega} \frac{Y^TP_X Y}{\sigma^2}\}]d \omega\\
\le &P(\omega \le \frac{t_n}{n}|Y)- \int_0^{\frac{t_n}{n}}\pi(\omega|\sigma^2,Y) \frac{(1-\omega)^{\frac{p}{4}}}{(1-\omega/2)^{\frac{p}{2}}} \exp\{ - \frac{1}{4}\cdot \frac{\omega^2}{2-\omega} \frac{Y^TP_X Y}{\sigma^2}\}d \omega.
\end{split}
\end{displaymath}
Here we let $\xi_n=\log n$ in $t_n = \xi_n p \log n$ . Then if $\omega \in (0, \frac{t_n}{n})$,
\begin{displaymath}
(1-\frac{t_n}{n})^{\frac{p}{4}} \exp\{ - \frac{1}{4}\cdot \frac{(\frac{t_n}{n})^2}{2-\frac{t_n}{n}} \frac{Y^TP_X Y}{\sigma^2}\} \le \frac{(1-\omega)^{\frac{p}{4}}}{(1-\omega/2)^{\frac{p}{2}}} \exp\{ - \frac{1}{4}\cdot \frac{\omega^2}{2-\omega} \frac{Y^TP_X Y}{\sigma^2}\} \le \frac{1}{(1-\frac{t_n}{2 n})^{\frac{p}{2}}}.
\end{displaymath}
Hence
\begin{displaymath}
\begin{split}
&(1-\frac{t_n}{n})^{\frac{p}{4}} \exp\{ - \frac{1}{4}\cdot \frac{(\frac{t_n}{n})^2}{2-\frac{t_n}{n}} \frac{Y^TP_X Y}{\sigma^2}\} P(\omega>\frac{t_n}{n}|Y)\\
\le & \int_0^{\frac{t_n}{n}}\pi(\omega|\sigma^2,Y) \frac{(1-\omega)^{\frac{p}{4}}}{(1-\omega/2)^{\frac{p}{2}}} \exp\{ - \frac{1}{4}\cdot \frac{\omega^2}{2-\omega} \frac{Y^TP_X Y}{\sigma^2}\}d \omega\\
\le & \frac{1}{(1-\frac{t_n}{2 n})^{\frac{p}{2}}} P(\omega>\frac{t_n}{n}|Y).
\end{split}
\end{displaymath}
By Lemma \ref{lemma3}, to show $E_0 I \to 0$,  it suffices to show 
\begin{displaymath}
P(\omega \le \frac{t_n}{n}|Y) E_0 [\exp\{ - \frac{1}{4}\cdot \frac{(\frac{t_n}{n})^2}{2-\frac{t_n}{n}} \frac{Y^TP_X Y}{\sigma^2}\}] \ge 1
\end{displaymath}  
as $n\to \infty$.
\begin{displaymath}
\begin{split}
&E_0 [P(\omega \le \frac{t_n}{n}|Y) \exp\{ - \frac{1}{4}\cdot \frac{(\frac{t_n}{n})^2}{2-\frac{t_n}{n}} \frac{Y^TP_X Y}{\sigma^2}\}]\\
\ge & E_0 [P(\omega \le \frac{t_n}{n}|Y) \exp\{ - \frac{1}{4}\cdot \frac{(\frac{t_n}{n})^2}{2-\frac{t_n}{n}} \frac{Y^TP_X Y}{\sigma^2}\}\textbf{1}_{S_n}]\\
\ge & E_0 [P(\omega \le \frac{t_n}{n}|Y) \exp\{ - \frac{1}{4}\cdot \frac{(\frac{t_n}{n})^2}{2-\frac{t_n}{n}} l_2 n\}\textbf{1}_{S_n}]\\
=& \exp\{ - \frac{1}{4}\cdot \frac{(\frac{t_n}{n})^2}{2-\frac{t_n}{n}} l_2 n\}E_0 [P(\omega \le \frac{t_n}{n}|Y)\textbf{1}_{S_n}],\\
\end{split}
\end{displaymath}
which implies 
\begin{displaymath}
\lim_{n\to \infty} E_0 [P(\omega \le \frac{t_n}{n}|Y) \exp\{ - \frac{1}{4}\cdot \frac{(\frac{t_n}{n})^2}{2-\frac{t_n}{n}} \frac{Y^TP_X Y}{\sigma^2}\}]\ge  \lim_{n\to \infty}E_0 [P(\omega \le \frac{t_n}{n}|Y)\textbf{1}_{S_n}] .
\end{displaymath}
We have already shown  $E_0 [P(\omega > \frac{t_n}{n}|Y)\textbf{1}_{S_n}] \to 0$  and we also have $E_0 \textbf{1}_{S_n} = P(S_n) \to 1$, as $n\to \infty$. So
\begin{displaymath}
E_0 [P(\omega \le \frac{t_n}{n}|Y)\textbf{1}_{S_n}]  \to 1,
\end{displaymath}
which  proves this theorem.\\
\\
\\
Let
\begin{displaymath}
\begin{split}
S_{n,1}&= \{Y: l_1 n \le   \frac{Y^TP_X Y}{\sigma_0^2}  \le l_2 n\} ,  \\
S_{n,2}&= \{Y: (n-p)-\sqrt{n-p} \log (n-p) \le \frac{Y^TP_X Y}{\sigma_0^2} \le (n-p)+\sqrt{n-p} \log (n-p)\},\\
\tilde{S}_n &= S_{n,1} \cap S_{n,2}.
\end{split}
\end{displaymath}
Then by Lemma \ref{lemma2}, we have $P(\tilde{S}_n )\to 1$, as $n\to \infty$.\\
\\
\textbf{Proof of lemma \ref{lemma4}:}\\
Similar to Lemma \ref{lemma3}, we just need to show 
\begin{displaymath}
\frac{\int_0^{\frac{t_n}{n}}\pi(\omega|Y)d \omega}{\int^1_{\frac{t_n}{n}}\pi(\omega|Y)d \omega} \textbf{1}_{\tilde{S}_n}  \to \infty  ,
\end{displaymath}
as $n\to \infty$.  In the following proof, we will use the notation $C(>0)$ as a generic constant.\\

\begin{equation}
\begin{split}
& \textbf{1}_{\tilde{S}_n}{\int_0^{\frac{t_n}{n}}(d+Y^T(I-P_X) Y+\omega Y^T P_X Y)^{-\frac{n+c}{2}}\omega^{a+\frac{p}{2}-1}(1- \omega)^{b-1} d \omega} \\
\ge & C {\int_0^{\frac{t_n}{n}}(d+\sigma^2_0((n-p)+\sqrt{n-p} \log (n-p))+\omega l_2 n)^{-\frac{n+c}{2}}\omega^{a+\frac{p}{2}-1} d \omega}  \\
\ge & C {\int_0^{\frac{\sqrt{t_n}}{n}}(d+\sigma^2_0((n-p)+\sqrt{n-p} \log (n-p))+\omega l_2 n)^{-\frac{n+c}{2}}\omega^{a+\frac{p}{2}-1} d \omega}  \\
\ge & C {(d+\sigma^2_0((n-p)+\sqrt{n-p} \log (n-p))+ l_2 \sqrt{t_n})^{-\frac{n+c}{2}}\int_0^{\frac{\sqrt{t_n}}{n}}\omega^{a+\frac{p}{2}-1} d \omega}  \\ 
=&C {(d+\sigma^2_0((n-p)+\sqrt{n-p} \log (n-p))+ l_2 \sqrt{t_n})^{-\frac{n+c}{2}}} \frac{(\frac{\sqrt{t_n}}{n})^{a+p/2}}{a+p/2}
\end{split} \label{numerator2}
\end{equation}
We also have
\begin{equation}
\begin{split}
&\textbf{1}_{\tilde{S}_n}{\int^1_{\frac{t_n}{n}}(d+Y^T(I-P_X) Y+\omega Y^T P_X Y)^{-\frac{n+c}{2}}\omega^{a+\frac{p}{2}-1}(1- \omega)^{b-1} d \omega} \\
\le & {(d+\sigma^2_0((n-p)-\sqrt{n-p} \log (n-p))+ l_1 {t_n})^{-\frac{n+c}{2}}} {\int^1_{\frac{t_n}{n}}\omega^{a+\frac{p}{2}-1}(1- \omega)^{b-1} d \omega} \\
\le & {(d+\sigma^2_0((n-p)-\sqrt{n-p} \log (n-p))+ l_1 {t_n})^{-\frac{n+c}{2}}} {\int^1_{0}\omega^{a+\frac{p}{2}-1}(1- \omega)^{b-1} d \omega} \\
= & {(d+\sigma^2_0((n-p)-\sqrt{n-p} \log (n-p))+ l_1 {t_n})^{-\frac{n+c}{2}}} \frac{\Gamma(b)\Gamma(a+p/2)}{\Gamma(a+b+p/2)}.
\end{split}
\label{denom 2}
\end{equation}
Then
\begin{displaymath}
\begin{split}
&\frac{\int_0^{\frac{t_n}{n}}\pi(\omega|Y)d \omega}{\int^1_{\frac{t_n}{n}}\pi(\omega|Y)d \omega} \textbf{1}_{\tilde{S}_n} \\
\ge & C (\frac{d+\sigma^2_0((n-p)-\sqrt{n-p} \log (n-p))+ l_1 {t_n}}{d+\sigma^2_0((n-p)+\sqrt{n-p} \log (n-p))+ l_2 \sqrt{t_n}})^{\frac{n+c}{2}}\frac{\Gamma(a+b+p/2)}{\Gamma(b)\Gamma(a+p/2)} \frac{(\frac{\sqrt{t_n}}{n})^{a+p/2}}{a+p/2} \\
\ge & C (1+\frac{-2\sigma^2_0 \sqrt{n-p} \log (n-p)+ l_1 {t_n}-l_2 \sqrt{t_n}}{d+\sigma^2_0((n-p)+\sqrt{n-p} \log (n-p))+ l_2 \sqrt{t_n}})^{\frac{n+c}{2}}\frac{\Gamma(a+b+p/2)}{\Gamma(b)\Gamma(a+p/2)} \frac{(\frac{\sqrt{t_n}}{n})^{a+p/2}}{a+p/2} \\
=& C (1+\frac{2}{n}\frac{-2\sigma^2_0 \sqrt{n-p} \log (n-p)+ l_1 {t_n}-l_2 \sqrt{t_n}}{\sigma^2_0/2+ \frac{l_2 \sqrt{t_n}+d-p\sigma^2_0+ \sqrt{n-p} \log (n-p) }{2n}})^{\frac{n+c}{2}}\frac{\Gamma(a+b+p/2)}{\Gamma(b)\Gamma(a+p/2)} \frac{(\frac{\sqrt{t_n}}{n})^{a+p/2}}{a+p/2}\\
\ge & C \exp(l_1 t_n \sigma_0^2)  \frac{\Gamma(a+b+p/2)}{\Gamma(b)\Gamma(a+p/2)} \frac{(\frac{\sqrt{t_n}}{n})^{a+p/2}}{a+p/2} \quad(assume \quad t_n \succeq \sqrt{n} \log n) \\
\to  &  \infty,
\end{split}
\end{displaymath}
as $n\to \infty$.\\
\\
\\
\textbf{Proof of Theorem \ref{theorem2}:}\\
 $\pi(\beta|Y)=\int_0^1\pi(\beta|\omega, Y)\pi(\omega|Y)d\omega$.\\
 
By convexity of $H^2$ and the Jensen inequality, we have 
\begin{displaymath}
\begin{split}
& H^2(\pi(\beta|Y),N(\hat{\beta}, \sigma_0^2(X^T X)^{-1}))\\
\le & \int_0^1\pi(\omega|Y) H^2(t_{n+c}((1-\omega)\hat{\beta},[\frac{(1-\omega)(n+c)}{d+Y^T(I-P_X) Y+\omega Y^T P_X Y}]^{-1}(X^TX)^{-1}) ,N(\hat{\beta},\sigma_0^2(X^T X)^{-1}))d \omega\\
\le &\int_0^{\frac{t_n}{n}}\pi(\omega|Y) H^2(t_{n+c}((1-\omega)\hat{\beta},[\frac{(1-\omega)(n+c)}{d+Y^T(I-P_X) Y+\omega Y^T P_X Y}]^{-1}(X^TX)^{-1}) ,N(\hat{\beta},\sigma_0^2(X^T X)^{-1}))   d \omega\\
& +\int_{\frac{t_n}{n}}^1\pi(\omega|Y) d \omega\\
= &\int_0^{\frac{t_n}{n}}\pi(\omega|Y) H^2(t_{n+c}((1-\omega)\hat{\beta},[\frac{(1-\omega)(n+c)}{d+Y^T(I-P_X) Y+\omega Y^T P_X Y}]^{-1}(X^TX)^{-1}) ,N(\hat{\beta},\sigma_0^2(X^T X)^{-1})) d \omega\\
& +P(\omega>\frac{t_n}{n}|Y)\\
= &III+IV(say).
\end{split}
\end{displaymath}
By Lemma \ref{lemma4}, we have $E_0 IV \to 0$, so it suffices to show $E_0 III \to 0$. \\
\\
To this end, we begin with
\begin{displaymath}
\begin{split}
& H^2(t_{n+c}((1-\omega)\hat{\beta},[\frac{(1-\omega)(n+c)}{d+Y^T(I-P_X) Y+\omega Y^T P_X Y}]^{-1}(X^TX)^{-1}) ,N(\hat{\beta},\sigma_0^2(X^T X)^{-1}))\\
\le & 2 H^2(t_{n+c}((1-\omega)\hat{\beta},[\frac{(1-\omega)(n+c)}{d+Y^T(I-P_X) Y+\omega Y^T P_X Y}]^{-1}(X^TX)^{-1}) ,N((1-\omega)\hat{\beta},\sigma_0^2(1-\omega)(X^T X)^{-1}))\\
&+ 2 H^2(N((1-\omega)\hat{\beta},\sigma_0^2(1-\omega)(X^T X)^{-1}) ,N(\hat{\beta},\sigma_0^2(X^T X)^{-1}))
\end{split}
\end{displaymath}
In Theorem \ref{theorem1}, we already showed that 
\begin{displaymath}
E_0 \int_0^{\frac{t_n}{n}}H^2(N((1-\omega)\hat{\beta},\sigma_0^2(1-\omega)(X^T X)^{-1}) ,N(\hat{\beta},\sigma_0^2(X^T X)^{-1})) \pi(\omega| Y) d \omega \to 0
\end{displaymath}
as $n \to \infty$. So we just need to show
\begin{scriptsize}
\begin{displaymath}
E_0 \int_0^{\frac{t_n}{n}}H^2(t_{n+c}((1-\omega)\hat{\beta},[\frac{(1-\omega)(n+c)}{d+Y^T(I-P_X) Y+\omega Y^T P_X Y}]^{-1}(X^TX)^{-1}) ,N((1-\omega)\hat{\beta},\sigma_0^2(1-\omega)(X^T X)^{-1})) \pi(\omega| Y) d \omega \to 0
\end{displaymath}
\end{scriptsize}
as $n \to \infty$. It suffices to show 
\begin{scriptsize}
\begin{displaymath}
E_0 \int_0^{\frac{t_n}{n}}[1-H^2(t_{n+c}((1-\omega)\hat{\beta},[\frac{(1-\omega)(n+c)}{d+Y^T(I-P_X) Y+\omega Y^T P_X Y}]^{-1}(X^TX)^{-1}) ,N((1-\omega)\hat{\beta},\sigma_0^2(1-\omega)(X^T X)^{-1})) ]\pi(\omega| Y) d \omega \to 1,
\end{displaymath}
\end{scriptsize}
as $n \to \infty$.\\
Here we let $\xi_n =\log n$ in $t_n = \xi_n p \log n$. When $\omega \in (0, \frac{t_n}{n})$, $Y \in \tilde{S}_n$,
\begin{displaymath}
\begin{split}
&  \frac{d+\sigma^2_0((n-p)-\sqrt{n-p} \log (n-p)+l_1 {n} \omega)}{n+c} \\
\le &   \frac{d+Y^T(I-P_X) Y+\omega Y^T P_X Y}{n+c}\\
\le &  \frac{d+\sigma^2_0((n-p)+\sqrt{n-p} \log (n-p)+l_2 {n} \omega)}{n+c}.
\end{split}
\end{displaymath}
By $1+x \le \exp(x)$, we also have
\begin{scriptsize}
\begin{displaymath}
\begin{split}
&  1-H^2(t_{n+c}((1-\omega)\hat{\beta},[\frac{(1-\omega)(n+c)}{d+Y^T(I-P_X) Y+\omega Y^T P_X Y}]^{-1}(X^TX)^{-1}) ,N((1-\omega)\hat{\beta},\sigma_0^2(1-\omega)(X^T X)^{-1}))  \\
=&  \frac{\Gamma^{1/2}(\frac{n+p+c}{2})}{\Gamma^{1/2}(\frac{n+p}{2})(\frac{n+c}{2})^{\frac{p}{4}}}(2 \pi)^{-\frac{p}{2}} |X^T X|^{\frac{1}{2}}(\sigma_0^2)^{-\frac{p}{4}}(1-\omega)^{-\frac{p}{2}}(\frac{d+Y^T(I-P_X) Y+\omega Y^T P_X Y}{n+c})^{\frac{p}{4}}\\
&\cdot \int  \exp(-\frac{(\beta-(1-\omega)\hat{\beta})^T(X^T X)(\beta-(1-\omega)\hat{\beta})}{4(1-\omega)\sigma_0^2})\\
&\cdot[1+\frac{1}{n+c}(\beta-(1-\omega)\hat{\beta})^T [\frac{(d+Y^T(I-P_X) Y+\omega Y^T P_X Y)(X^T X)}{(1-\omega)(n+c)}]^{-1}(\beta-(1-\omega)\hat{\beta})]^{-\frac{n+p+c}{4}}  d \beta \\
\succeq &   \frac{\Gamma^{1/2}(\frac{n+p+c}{2})}{\Gamma^{1/2}(\frac{n+p}{2})(\frac{n+c}{2})^{\frac{p}{4}}}(2 \pi)^{-\frac{p}{2}} |X^T X|^{\frac{1}{2}}(\sigma_0^2)^{-\frac{p}{4}}(1-\omega)^{-\frac{p}{2}}\\
& \cdot  \int  \exp(-\frac{(\beta-(1-\omega)\hat{\beta})^T(X^T X)(\beta-(1-\omega)\hat{\beta})}{4(1-\omega)\sigma_0^2}) \\
& \cdot    \exp(-\frac{(\beta-(1-\omega)\hat{\beta})^T(X^T X)(\beta-(1-\omega)\hat{\beta})}{4(1-\omega)\sigma_0^2} \frac{n+p+c}{n+c}) \frac{(n-p)+\sqrt{n-p}+d+l_2 n \omega}{n+c}   d \beta  \\
= & \frac{\Gamma^{1/2}(\frac{n+p+c}{2})}{\Gamma^{1/2}(\frac{n+p}{2})(\frac{n+c}{2})^{\frac{p}{4}}}(\frac{1}{2}+\frac{n+p+c}{2(n+c)}\frac{(n-p)+\sqrt{n-p}+d+l_2 t_n }{n+c})^{-\frac{p}{2}}.
\end{split}
\end{displaymath}
\end{scriptsize}
Then we have
\begin{scriptsize}
\begin{displaymath}
\begin{split}
&  E_0 \int_0^{\frac{t_n}{n}}[1-H^2(t_{n+c}((1-\omega)\hat{\beta},[\frac{(1-\omega)(n+c)}{d+Y^T(I-P_X) Y+\omega Y^T P_X Y}]^{-1}(X^TX)^{-1}) ,N((1-\omega)\hat{\beta},\sigma_0^2(1-\omega)(X^T X)^{-1})) ]\pi(\omega| Y) d \omega \\
\ge & E_0 \int_0^{\frac{t_n}{n}}[1-H^2(t_{n+c}((1-\omega)\hat{\beta},[\frac{(1-\omega)(n+c)}{d+Y^T(I-P_X) Y+\omega Y^T P_X Y}]^{-1}(X^TX)^{-1}) ,N((1-\omega)\hat{\beta},\sigma_0^2(1-\omega)(X^T X)^{-1})) ]\pi(\omega| Y) \textbf{1}_{\tilde{S}_n} d \omega \\
\ge & \frac{\Gamma^{1/2}(\frac{n+p+c}{2})}{\Gamma^{1/2}(\frac{n+p}{2})(\frac{n+c}{2})^{\frac{p}{4}}}(\frac{1}{2}+\frac{n+p+c}{2(n+c)}\frac{(n-p)+\sqrt{n-p}+d+l_2 t_n }{n+c})^{-\frac{p}{2}} E_0 \int_0^{\frac{t_n}{n}}\pi(\omega| Y) \textbf{1}_{\tilde{S}_n} d \omega  \to 1
\end{split}
\end{displaymath}
\end{scriptsize}
as $n \to \infty$. This proves Theorem 2.\\
\\
\noindent
\textbf{Proof of Theorem \ref{theorem3}:}\\
By Lemma 4 in supplementary material of \citep{Johnson:2012}, when the conditions of Theorem 1 in \citep{Johnson:2012} hold, we have
\begin{equation}
Q_t=E_t(\prod_{i=1}^t \beta^{2r}_{t_i}) \stackrel{a.s}{\rightarrow} \prod_{i=1}^t (\beta^0_{t_i})^{2r}.
\end{equation}
Thus we have
\begin{equation}
 H^2(\pi(\beta_t |\tau, r, \sigma^2, Y),N(\tilde{\beta}_t, \sigma^2(C_t)^{-1})) \to 0 ,\label{begin}
\end{equation}
by showing 
\begin{displaymath}
\int \sqrt{\pi(\beta_t |\tau, r, \sigma^2, Y)N(\beta_t| \tilde{\beta}_t, \sigma^2(C_t)^{-1})} d \beta_t \to 1,
\end{displaymath}
as $n \to \infty$. \\

We begin with 
\begin{equation}
H^2(N(\hat{\beta}_t, \sigma^2(C_t)^{-1}),N(\tilde{\beta}_t, \sigma^2(C_t)^{-1}))=1-\exp\{-\frac{1}{8}(\tilde{\beta}_t-\hat{\beta}_t)^T \sigma^{-2} C_t (\tilde{\beta}_t-\hat{\beta}_t)\}.
\end{equation}
But
\begin{displaymath}
0<C \le \lambda_{min}(\frac{X_t^T X_t}{n}) \le \lambda_{max}(\frac{X_t^T X_t}{n}) \le M < \infty 
\end{displaymath}
and
\begin{displaymath}
0<a_1 \le \lambda_{min}(A_t) \le \lambda_{max}(A_t) \le a_2 < \infty,
\end{displaymath}
also
\begin{displaymath}
\tilde{\beta}_t-\hat{\beta}_t =-(X_t^T X_t)^{-1}(\tau A_t^{-1}+(X_t^T X_t)^{-1})^{-1}(X_t^T X_t)^{-1} X_t^T Y.
\end{displaymath}
Hence, plug it in the following equation we get
\begin{displaymath}
\begin{split}
& (\tilde{\beta}_t-\hat{\beta}_t)^T  C_t (\tilde{\beta}_t-\hat{\beta}_t) \\
\le & (\frac{1}{Cn})^2 Y^T X_t(X_t^T X_t)^{-1} (\tau A_t^{-1}+(X_t^T X_t)^{-1})^{-2} (X_t^T X_t)^{-1}X_t^T Y\\
\le & (\frac{1}{Cn})^2 (\frac{1}{Cn}) (\frac{a_2}{\tau})^2(\frac{1}{Cn}) Y^T_n X_t X_t^T  Y_n \\
\preceq & tr(X_t X_t^T)\frac{Y^T_n Y_n}{n^4} \preceq  \frac{Y^T_n Y_n}{n^3} \stackrel{p}{\rightarrow}0.
\end{split}
\end{displaymath}
Thus,
\begin{displaymath}
H^2(N(\hat{\beta}_t, \sigma^2(C_t)^{-1}),N(\hat{\beta}_t, \sigma^2(C_t)^{-1})) \stackrel{p}{\rightarrow}0 .
\end{displaymath}
Now we consider 
\begin{equation}
H^2(N(\hat{\beta}_t, \sigma^2(C_t)^{-1}),N(\tilde{\beta}_t, \sigma^2(X_t^T X_t)^{-1}))=1-\frac{1}{det^{\frac{1}{4}}(C_t \frac{C_t^{-1}+(X_t^T X_t)^{-1}}{2})det^{\frac{1}{4}}((X_t^T X_t)\frac{C_t^{-1}+(X_t^T X_t)^{-1}}{2})},
\end{equation}
since 
\begin{displaymath}
|C_t \frac{C_t^{-1}+(X_t^T X_t)^{-1}}{2}|=|I_t + \frac{A_t}{2}(X_t^T X_t)^{-1}|,
\end{displaymath}
and 
\begin{displaymath}
1 \le |I_t + \frac{A_t}{2}(X_t^T X_t)^{-1}| \le  |I_t + \frac{a_2}{2nc}I_t| =(1+\frac{a_2}{2nc})^t \to 1,
\end{displaymath}
we get
\begin{displaymath}
|C_t \frac{C_t^{-1}+(X_t^T X_t)^{-1}}{2}| \to 1.
\end{displaymath}
Similarly, we have 
\begin{displaymath}
|(X_t^T X_t)\frac{C_t^{-1}+(X_t^T X_t)^{-1}}{2}| \to 1.
\end{displaymath}
Thus 
\begin{equation}
H^2(N(\hat{\beta}_t, \sigma^2(C_t)^{-1}),N(\tilde{\beta}_t, \sigma^2(X_t^T X_t)^{-1})) \stackrel{p}{\rightarrow}0. \label{end}
\end{equation}
By (\ref{begin})-(\ref{end}), we have 
\begin{equation}
 H^2(\pi(\beta_t |\tau, r, \sigma^2, Y),N(\hat{\beta}_t, \sigma^2(X^T_t X_t)^{-1})) \stackrel{p}{\rightarrow}0 .
\end{equation}
Since 
\begin{displaymath}
\pi(\beta|\tau, r, \sigma^2, Y)=\sum_k \pi(\beta_k|k,\tau, r, \sigma^2, Y) \pi(k|Y),
\end{displaymath}
by convexity of $H^2$ and the fact that $\pi(t|Y) \stackrel{p}{\rightarrow}1 $  \citep{Johnson:2012} Theorem 1, we have
\begin{displaymath}
E_{\beta_t^0} H^2(\pi(\beta|\tau, r, \sigma^2, Y),N(\hat{\beta}_t, \sigma^2(X_t^T X_t)^{-1})) \to 0 ,
\end{displaymath}
as $n \to \infty$.

\bibliography{sample}

\end{document}